\documentclass[12pt]{article}
\usepackage{amsmath,amsfonts,amsthm,amssymb}
\def\N{{\mathbb N}}
\def\R{{\mathbb R}}
\def\Z{{\mathbb Z}}
\def\Ce{{\mathbb C}}

\newtheorem{theorem}{Theorem}

\newtheorem{lemma}[theorem]{Lemma}

\def\d{\partial}

\begin{document}

\title{Spectral projections for the twisted Laplacian}

\author{Herbert Koch
    \\ Fachbereich Mathematik   \\ Universit\"at Dortmund
\and Fulvio Ricci \\
    Scuola Normale Superiore\\Pisa}

\maketitle

\begin{abstract} Let $n \ge 1$, $d=2n$, $(x,y)\in \R^n\times \R^n$
the notation for a generic
point in $\R^{2n}$. The twisted Laplacian
\[
L = -\frac12 \sum_{j=1}^n \left[ (\d_{x_j} + iy_j)^2 + (\d_{y_j} - i x_j)^2
\right]
    \]
has the  spectrum
\[
\{ n + 2k =\lambda^2: k \text{ a nonnegative integer }\}.
\] Let
$\mu=\lambda^2$ and
let $P_\mu$ be the spectral projection onto the (infinite dimensional)
eigenspace. We find the
optimal exponent $\rho(p)$ in the estimate
\[ \Vert P_\lambda u \Vert_{L^p(\R^d)} \lesssim  \lambda^{\rho(p)}
\Vert  u \Vert_{L^2(\R^d)}
\]
for all $p\in[2,\infty]$, improving previous partial results by
Thangavelu and by Stempak
and Zienkiewicz.  The expression for $\rho(p)$ is
\[
\rho(p) = \left\{ \begin{array}{cl}
\frac1p - \frac12 & \text{ if }  2 \le p \le \frac{2(d+1)}{d-1} \\
\frac{d-2}2 - \frac{d}p & \text{ if } \frac{2(d+1)}{d-1} \le p \le \infty\ .
\end{array}
\right.
\]
\end{abstract}

\section{Introduction} Sharp $L^2-L^p$-bounds for spectral
projections onto eigenspaces of
elliptic differential operators $L$ with a discrete spectrum have
attracted considerable
attention in the last 20 years, starting with the work of C. Sogge
\cite{MR94c:35178} on the spherical
Laplacian (see also \cite{MR89d:35131}). If $P_\lambda$ is the spectral
projection corresponding to an
eigenvalue $\lambda^2$ of $L$ and $2\le p\le\infty$, one looks for
the best possible
exponent $\rho(p)$ such that
\begin{equation}\label{exponent}
\|P_\lambda
u\|_p\lesssim \lambda^{\rho(p)}\|u\|_2\ .
\end{equation}

In general, $\rho(p)$ is a convex function of
$1/p$. Strict convexity of $\rho(p)$ at some points is connected with
dispersive
estimates for $L$ and in
some cases to phenomena in harmonic analysis such as restriction
theorems for the Fourier
transform.

We consider
here the so-called twisted Laplacian
\[ L = -\frac12 \sum_{j=1}^n \left[ (\d_{x_j} + iy_j)^2 + (\d_{y_j} - i x_j)^2
\right]
    \]
in $\R^{2n}$. In dimension 2 (i.e. $n=1$), $L$ can be viewed as a
Schr\"odinger operator with
the magnetic potential $A=(y,-x)$, inducing a constant magnetic field.

The sharp estimates \eqref{exponent} relative to $L$ are as follows.

\begin{theorem}
With $d=2n$ and $2\le p\le\infty$, \eqref{exponent} holds with
\begin{equation} \label{rho}
\rho(p) = \left\{ \begin{array}{cl}
\frac1p - \frac12 & \text{ if }  2 \le p \le \frac{2(d+1)}{d-1} \\
\frac{d-2}2 - \frac{d}p & \text{ if } \frac{2(d+1)}{d-1} \le p \le \infty\ ,
\end{array}
\right.
\end{equation}
and with no smaller exponent.
\end{theorem}

The first partial result in this direction is in \cite{thangavelu},
where the second exponent in
\eqref{rho} was obtained for $p$ larger than some $p_0>\frac{2(d+1)}{d-1}$.
Later,  Stempak and Zienkiewicz
\cite{MR2000e:35166} obtained \eqref{rho} for $p\ne\frac{2(d+1)}{d-1}$.

Our approach is inspired by the recent work of Koch and Tataru
\cite{oi} on dispersive estimates and their application to
the Hermite operator \cite{hermite}. In particular it based on PDE techniques
and it does  not need the refined estimates on
Laguerre functions, which are used in \cite{MR2000e:35166}.

The basic estimate is the local dispersive estimate
\eqref{dispersive2} below, which implies
the endpoint result for $p=\frac{2(d+1)}{d-1}$ by a simple covering
argument. Weaker local
estimates were also used by Thangavelu. The other endpoint result at
$p=\infty$ follows
from the exact computation of the $L^2-L^\infty$-norm of $P_\lambda$,
see \eqref{twoinfty}
below.

The connection between $L$ and Hermite operators is two-fold. On one
hand, if we decompose $L^2(\R^{2n})$ as the orthogonal sum of
the subspaces $V_m$
consisting of the functions $f$ such that
$$
f\big(e^{i\theta_1}(x_1+iy_1),\dots,e^{i\theta_n}(x_n+iy_n)\big)=e^{i(m_1\theta_1+\cdots+m_n
\theta_n)}f(x_1+iy_1,\dots,x_n+iy_n)\ ,
$$
($m\in\Z^n$), we see that
\[ L = -\frac12\big( \Delta_{\R^{2n}} -|x|^2-|y|^2\big)+i\sum_{j=1}^n
(x_j\d_{y_j}-y_j\d_{x_j}  )\ ,
    \]
equals $-\frac12\big( \Delta_{\R^{2n}} -|x|^2-|y|^2\big)-\sum m_j$ on $V_m$.

On the other hand, the operators $X_j=\d_{x_j}+iy_j$ and
$Y_j=\d_{y_j}-ix_j$ satisfy the
canonical commutation relations $[X_j,Y_k]=-2i\delta_{j,k}$. This
implies that there is a
unitary projective representation $\pi$ of $\R^{2n}$ on $L^2(\R^n)$
(with variable $\xi$),
called the Weyl representation, and such that $d\pi(X_j)=\sqrt2\d_{\xi_j}$,
$d\pi(Y_j)=-i\sqrt2\xi_j$. Then $d\pi(L)=-\Delta_{\R^n} +|\xi|^2$.
The Stone-von~Neumann theorem establishes an isomorphism between
$L^2(\R^{2n})$ and
$L^2(\R^n)\otimes L^2(\R^n)$, which intertwines the action of the
twisted Laplacian on the
first space with the action of the Hermite operator on the first
factor in the tensor product.

Each of these arguments shows that $L$ has a discrete spectrum, equal
to the set of integers
$\lambda^2=n+2k$, $k\in\N$, and that each eigenspace is
infinite-dimensional. The precise
description of the eigenspaces is given in Section~2.

The twisted Laplacian also describes the action of the Heisenberg
sub-Laplacian on special
classes of functions. On the reduced Heisenberg group
$h_n=\R^n\times\R^n\times\mathbb{T}$
with product
$$
(x,y,e^{i\theta})(x',y',e^{i\theta'})=(x+x',y+y',e^{i(\theta+\theta'+x'\cdot
y-x\cdot y')})\ ,
$$
consider the left-invariant sub-Laplacian
$$
{\cal
L}=-\frac12\sum_{j=1}^n\left[(\d_{x_j}+y_j\d_\theta)^2+(\d_{y_j}-x_j\d_\theta)^2
\right]\ .
$$

If $f(x,y,e^{i\theta})=g(x,y)e^{im\theta}$, with $m\in\Z$, then
\begin{eqnarray*}
{\cal L}f(x,y,e^{i\theta})&=& -\frac12e^{im\theta} \sum_{j=1}^n
\left[ (\d_{x_j} + imy_j)^2 +
(\d_{y_j} - i mx_j)^2
\right]g(x,y)\\
&=&e^{im\theta}L_mg(x,y)\ .
\end{eqnarray*}

One easily verifies that, for $m\ne0$, the spectrum of $L_m$ consists
of the integers
$|m|(n+2k)$, with $k\in\N$, and that the pairs $(|m|(n+2k),m)$, with
$m\in\Z\setminus\{0\},k\in\N$, give the discrete joint spectrum of
$\cal L$ and
$i^{-1}\d_\theta$.
If $P_{m,k}$ is the orthogonal projection on the joint eigenspace, a
simple scaling argument
shows that
$$
\|P_{m,k}u\|_{L^p(h_n)}\lesssim (n+2k)^{\frac12\rho(p)}|m|^{\sigma(p)}
\Vert u \Vert_{L^2(h_n)} \ ,
$$
with $\rho(p)$ as in \eqref{rho} and $\sigma(p)=\frac d2(\frac12-\frac1p)$.

{\bf Acknowledgments.}
This research  was done  during the Research Trimester on Harmonic Analysis at
the Centro de Giorgi in Pisa. Both authors are members of the European IHP
network  HARP ``Harmonic Analysis and Related Problems".
  We acknowledge the support for this work by
the European Commission
through the Network  and the support by the Centro de Giorgi in Pisa.

\section{Spectrum, eigenfunctions and lower bounds for the spectral
projections} Here we
introduce some notation, derive formulas for some eigenfunctions  and
calculate $L^p$ norms
of them. This will imply that the bounds are sharp.

We identify $\R^n \times \R^n$ with $\Ce^n$ by $z=x+iy$ and denote
\[ \d_{z_j} = \frac12( \d_{x_j} - i \d_{y_j}), \qquad \d_{\bar z_j} =
\frac12( \d_{x_j} + i
\d_{y_j}). \] Then
\[  -\frac12 \big[ (\d_{x_j} + i y_j)^2 + (\d_{y_j} - i x_j)^2 \big]
=  - \frac12 (2\d_{z_j}
-  \overline{z_j})(2\d_{\bar z_j} +  z_j)+1  =  \frac12 D_j^* D_j + 1 \ ,
\]
with $D_j=  2\d_{\bar z_j}+z_j$. Since
\[
    D_j (f(z,\bar z) e^{-\frac{|z|^2}2}) =  2(\d_{\bar z_j} f)
e^{-\frac{|z|^2}{2}},
\]
for all holomorphic function $f$
\[    L\left[f(z) e^{-\frac{|z|^2}2}\right] = nf(z)e^{-\frac{|z|^2}2} \ .
\] On the other hand
\[
\int_{\Ce^n} (L u) \bar  u dx = n \Vert u \Vert_{L^2}^2 + \frac12
\sum_{j=1}^n  \Vert D_j u \Vert_{L^2}^2\ ,
\]  which shows that $n$ is the smallest eigenvalue and that the
eigenfunctions  to that
eigenvalue are anniliated by $D_j$. An easy calculation shows
\[
    [D_j ,D_j^*] = 4
\] and hence, if $f$ is an eigenfunction to the eigenvalue $\mu$ then
\[  L D_j^* f = D_j^*  Lf + 2 D_j^* f = (\lambda + 2) D_j^* f \ ,
\] and
\[  L D_j f = D_j L f - 2 D_j f = (\lambda -2)D_j f\ .
\]

We obtain an orthogonal basis of
eigenfunctions of the form
\begin{equation}\label{eigen}  f_{\alpha,\beta} =(-2^{-1}D^*)^\alpha(z^\beta
e^{-\frac{|z|^2}2})= (-1)^{|\alpha|+|\beta|}e^{\frac{|z|^2}2}
\d_z^\alpha
\d_{\bar z}^\beta  e^{-|z|^2}
\end{equation}
with
\[   L f_{\alpha, \beta} = (n+2|\alpha|) f_{\alpha, \beta}\ .
\] In particular,
\[  \bar z_1^k e^{-|z|^2/2}= (-1)^kf_{(k,0,\dots),0}(z)
\] is an eigenfunction to the eigenvalue $n+2k$.  When we consider
the twisted Laplacian as a
quantization of the motion  of a charged particle in a constant
magnetic field, then this
eigenfunctions
    corresponds to the motion
    in a circle of radius $\sqrt{k}$ around zero and its maximal
concentration (to scale $1$)
    around that circle is dictated by the uncertainty principle.  Thus
this eigenfunction
behaves like a characteristic function of a neighborhood of size $1$
around that circle.

Its $L^p$ norm can be explicitly computed. Integrating in $z_1$
first, we obtain
\begin{eqnarray*}
\label{norm}
\int_{\Ce} |z_1|^{kp} e^{-\frac{p}2|z_1|^2} dz_1 &= & 2\pi
\int_0^\infty r^{kp+1}
e^{-\frac{p}2 r^2} dr \\ & = & \frac{2\pi}p
\left(\frac{2}p\right)^{\frac{kp}2} \int_0^\infty
t^{\frac{kp}2} e^{-t} dt \\ & = &  \frac{2\pi}p
\left(\frac{2}p\right)^{\frac{kp}2}
\Gamma(\frac{kp}2+1)
\end{eqnarray*}
and hence
\[
\Vert z_1^k e^{-\frac{|z|^2}2} \Vert_{L^p(\R^{2n})}^p  =  (\frac{2\pi}p)^n
\left(\frac{2}p\right)^{\frac{kp}2} \Gamma(\frac{kp}2+1)
\]

By Stirling's formula
\[
\Gamma(t+1) \sim \sqrt{2\pi t} (t/e)^t
\]
   and
\begin{equation} \label{pnorm}
\Vert z_1^k e^{-\frac{|z|^2}2} \Vert_{L^p(\R^{2n})} /  \Vert z_1^k
e^{-\frac{|z|^2}2}
\Vert_{L^2(\R^{2n})}
\sim  k^{\frac12(\frac1p - \frac12)}\ .
\end{equation}

This proves that $\rho(p)\ge \frac1p-\frac12$.

\vskip.3cm

We shall now need the $L^2$-norms of all the $f_{\alpha,\beta}$. We can reduce
ourselves to one dimension, since
$$
f_{\alpha,\beta}(z)=\prod_{j=1}^nf_{\alpha_j,\beta_j}(z_j)\ .
$$

In dimension $n=1$, the
computation in \eqref{pnorm} shows that
$$
\|f_{0,\ell}\|_2^2=\pi\ell!\ .
$$

Observe next that
\begin{eqnarray*}
\|f_{k,\ell}\|^2_2&=&\frac14\|D^*f_{k-1,\ell}\|_2^2\\
&=&\frac14\langle f_{k-1,\ell},DD^*f_{k-1,\ell}\rangle\\
&=&\frac14\langle f_{k-1,\ell},(D^*D+[D,D^*])f_{k-1,\ell}\rangle\\
&=&\frac14\langle f_{k-1,\ell},(2L+2)f_{k-1,\ell}\rangle\\
&=&k\|f_{k-1,\ell}\|_2^2\ ,
\end{eqnarray*}
so that $
\|f_{k,\ell}\|_2^2=\pi k!\ell!$,
and, in $n$ dimensions,
\begin{equation}\label{l2norms}
\|f_{\alpha,\beta}\|_2^2=\pi^n\alpha!\beta!\ .
\end{equation}

As for the Hermite operator \cite{hermite}, we expect that radial
eigenfunctions
will provide the sharp value of $\rho(p)$ for $p$ close to $\infty$.

For $k\in\N$, consider
\begin{eqnarray*}
f_k(z)&=& e^{\frac{|z|^2}2}\Delta^ke^{-|z|^2}\\
&=&4^{-k}e^{\frac{|z|^2}2}\Big(\sum_{j=1}^n\d_{z_j}\d_{\bar
z_j}\Big)^ke^{-|z|^2}\\
&=&4^{-k}\sum_{|\alpha|=k} \binom{k}{\alpha}f_{\alpha,\alpha}(z)\ .
\end{eqnarray*}

Then $f_k$ is radial, and it is an eigenfunction of $L$ with
eigenvalue $n+2k$. Since for
radial functions the twisted Laplacian and the Hermite operator
coincide up to a factor
$2$, the $f_k$ are the unique (up to scalar multiples) radial
eigenfunctions of $L$.

By \eqref{l2norms} and the orthogonality of the $f_{\alpha,\alpha}$,
\begin{eqnarray*}
\|f_k\|_2^2&=&4^{-2k} \pi^n\sum_{|\alpha|=k}k!^2\\
&=&4^{-2k} \pi^nk!^2\binom{n+k-1}{ k}\ .
\end{eqnarray*}

   From \eqref{eigen} we obtain that $f_{\alpha,\alpha}(0)=\alpha!$, hence
$$
f_k(0)= 4^{-k} k!\binom{n+k-1}{ k}\ .
$$

        Thus
\begin{eqnarray*}
\frac{\|f_k\|_{L^\infty}}{\Vert f_k \Vert_{L^2} }  & \ge &
\frac{f_k(0)}{\Vert f_k \Vert_{L^2} } \\ & = & \pi^{-\frac n2} \sqrt{
\binom{n+k-1}{k}}
\\ &\sim &
    k^{\frac{d-2}{4}}.
\end{eqnarray*}

This proves that $\rho(\infty)\ge \frac{d-2}2$.

It not hard to see that the first inequality is in fact an equality,
and that the ratio
$\|f_k\|_{L^\infty}/\| f_k\|_{L^2}$ coincides with the
$L^2-L^\infty$-norm of the spectral
projection. If
$f$ is an eigenfunction  to the eigenvalue
$n+2k$ of the twisted Laplacian with $L^2$ norm $1$, which  maximizes
the $L^\infty$ norm
then we may assume, after a twisted translation of the form
$f(x,y)\longmapsto e^{i(a\cdot
y-b\cdot x)}f(x-a,y-b)$, that it assumes its maximum at
$z=0$. Averaging over the unitary group
$U(n)$ we see that it has to be radial, hence a scalar multiple of
$f_k$. Thus, for
$\lambda^2=n+2k$,
\begin{equation}\label{twoinfty} \Vert P_\lambda \Vert_{2\to \infty}
=  \pi^{-n} \sqrt{\binom{n+k-1}{k}}\ .
\end{equation}

We look now for an estimate from below of $\|f_k\|_p$ for $p$ finite.
Since $f_k$ is radial we obtain by the divergence theorem, if $|z|=r$,
\begin{eqnarray*}
\d_rf_k(z)&=&\frac{c_d}{r^{d-1}}\int_{S_r}\d_\nu f_k(w)\,d\sigma(w)\\
&=&\frac{c_d}{r^{d-1}}\int_{B_r}\Delta f_k(w)\,dw\\
&=&\frac{c_d}{r^{d-1}}\int_{B_r}(-2L+|w|^2) f_k(w)\,dw\\
&=&\frac{c_d}{r^{d-1}}\int_{B_r}(-2n+4k+|w|^2) f_k(w)\,dw\\
\end{eqnarray*}

Thus, for $r<1$,
$$
|\d_rf(z)|\le Cr(n+2k)\|f_k\|_\infty=Cr(n+2k)f_k(0)\ ,
$$
and
$$
f(z)\ge f(0)\big(1-C(n+2k)|z|^2\big)\ .
$$

It follows that for $|z|<r_k=1/\sqrt{C(n+2k)}$, $f_k(z)>\frac12f(0)$, so that
$$
\|f_k\|_p\ge C  r_k^{\frac dp}f_k(0)\sim k^{-\frac
d{2p}+\frac{d-2}4}\|f_k\|_2\ .
$$

This shows that $\rho(p)\ge \frac{d-2}2-\frac dp$.

\section{The upper bounds for the spectral projections}

Let $\lambda^2 = n+2k$   and
\[ L u = \lambda^2 u. \]
Because of \eqref{twoinfty} the
   assertion follows once we prove for all
$(x_0,y_0)$ that
\begin{equation} \label{dispersive}
\lambda^{\frac1{d+1}}
\Vert u \Vert_{L^{\frac{2(d+1)}{d-1}}(B_{\lambda}(x_0,y_0))}
\lesssim
    \Vert u \Vert_{L^2(B_{2\lambda}(x_0,y_0))}.
\end{equation}

We set
\[ x =  x_0 + \lambda \bar x, \quad y = y_0 + \lambda \bar y \quad
    \bar u(\bar x, \bar y) = e^{i (x_0 y-y_0 x)}    u(x,y)
\] Then
\[   (\lambda^{-1}  \d_{\bar x_j} + i\lambda \bar y_j) \bar u  =
(\d_{x_j} + (y_j-(y_0)_j )
e^{i(x_0y-y_0 x)} u  =   e^{i(x_0y-y_0 x)} (\d_{x_j} + i y_j ) u
\] and hence
\[
\bar L \bar u := -\frac12 \sum_{j=1}^n
\left( (\d_{\bar x_j}+ i\lambda^2 \bar y_j)^2  + (\d_{\bar y_j}-
i\lambda^2 \bar x_j)^2
\right) \bar u =
\lambda^4 \bar u.
\] We drop the bar on $x$, $y$ and $u$ but not on $L$. Hence we study
\[  \bar L u = \lambda^4u \] in ball of radius $2$.  The inequality
\eqref{dispersive} takes the  form
\begin{equation}\label{dispersive2}
\lambda^{\frac2{d+1}} \Vert u \Vert_{L^{\frac{2(d+1)}{d-1}}(B_1(0))}
\lesssim
    \lambda  \Vert u \Vert_{L^2(B_2(0))}
\end{equation}

Actually a  slightly stronger bound than \eqref{dispersive2} is true:
\begin{lemma} \label{dispersivelem}   Suppose that
\[  L_{\mu} u := -\frac12 \sum_{j=1}^n
\left( (\d_{\bar x_j}+ i\mu \bar y_j)^2  + (\d_{\bar y_j}+ i\mu \bar
x_j)^2 \right) u - \mu^2
u = f \] Then
\[
\mu^{\frac1{d+1}}
    \Vert u \Vert_{L^{\frac{2(d+1)}{d-1}}(B_1(0))}
\lesssim   \mu^{1/2} \Vert u \Vert_{L^2(B_2(0))}  + \mu^{-1/2} \Vert
f \Vert_{L^2(B_2(0))}
\]
\end{lemma}
The dispersive estimate \eqref{dispersive2} is an immediate consequence.

\begin{proof}
The statement of Lemma \ref{dispersivelem} follows from Theorem 3 of
\cite{oi} in the same way as Lemma 3.4  of \cite{hermite} is deduced from the
same result.  The symbol of $L_\mu$ is
\[p(x,y,\xi,\eta)= \frac12( |\xi +y|^2+ |\eta-x|^2) - \mu^2 \]
which is real. Keeping $x$ and $y$ fixed it vanishes on a sphere of radius
$\mu$, which has $2n-1$ nonvanishing curvatures of size $\mu^{-1}$. This
curvature leads to the desired estimate.
\end{proof}


\end{document}